\documentclass[a4paper,11pt]{article}
\newtheorem{theo}{Theorem}[section]
\newtheorem{ex}[theo]{Example}
\newtheorem{prop}[theo]{Proposition}
\newtheorem{lem}[theo]{Lemma}
\newtheorem{cor}[theo]{Corollary}

\newtheorem{defi}[theo]{Definition}

\newtheorem{rema}[theo]{Remark}

\newtheorem{theobis}{Theorem}

\usepackage{amscd}
\usepackage{amssymb}
\usepackage{amsmath}
\usepackage{times}

\def\bb{\mathbb}

\title{ Kac-Moody algebras and Lie algebras of regular vector fields on tori}
\author{M. Rausch de Traubenberg
$\,\,$${}^{a,b}$ and
M.J. Slupinski
$\,\,$${}^{c}$\\
\\
{\small ${}^{a}${\it Laboratoire de Physique
Math\'ematique et Th\'eorique,
Universit\'e Montpellier II,}}\\
{\small {\it Place E. Bataillon,  34095 Montpellier,
France} }\\
{\small ${}^{b}${\it
Laboratoire de Physique Th\'eorique, 3 rue de
l'Universit\'e, 67084 Strasbourg, France}}\\
{\small  e-mail: rausch@lpm.univ-montp2.fr}\\
{\small ${}^c${\it Institut de Recherches en Math\'ematique
Avanc\'ee}}\\
{\small { \it Universit\'e Louis-Pasteur, and CNRS}}\\
{\small {\it 7 rue R. Descartes, 67084 Strasbourg Cedex, France}}\\
{\small  e-mail: slupins@math.u-strasbg.fr}}

\date{\today}

\begin{document}
\baselineskip=15pt
\maketitle

{\bf Abstract}: {\it We consider the problem of representing the Kac-Moody
algebra $\mathfrak{g}(N)$ specified by an $r\times
r$ indecomposable   generalised Cartan matrix $N$ as vector fields on the
torus ${{\bb C}^*}^r$. It is shown that,
if the  representations are of a certain form, this is possible if and only
if $\mathfrak{g}(N)\cong sl(r+1,{\bb C})$ or
$\tilde{sl}(r,{\bb C})$. For $sl(r+1,{\bb C})$ and $\tilde{sl}(r,{\bb C})$,
discrete families of
representations  are constructed. These generalise the well-known discrete
families of representations
of $sl(2,{\bb C})$ as regular vector fields on ${\bb C}^*$. }

\section{Introduction}
Let $H,X_1,X_{-1}$ be a standard basis of the Lie algebra  $sl(2,{\bb C})$.
As is well-known,
$sl(2,{\bb C})$ can be represented by regular vector fields on ${\bb C}^*$
according to the formulae:
$$
R(H)=2z\frac{d}{dz},\quad R(X_1)=z(z\frac{d}{dz}),\quad
R(X_{-1})=\frac{1}{z}(-z\frac{d}{dz}),
$$
or more generally $(n\in{\bb Z}^*)$ according to:
$$
R_n(H)=2\frac{z}{n}\frac{d}{dz},\quad
R_n(X_1)=z^n(\frac{z}{n}\frac{d}{dz}),\quad
R_n(X_{-1})=\frac{1}{z^{n}}(-\frac{z}{n}\frac{d}{dz}).
$$
Replacing $z$ by $e^{i\theta}$  gives representations of $sl(2,{\bb C})$ by
complex vector fields on the circle $U(1).$
 In this paper we consider the problem of obtaining analogous
families of representations by vector fields on ${{\bb C}^*}^r$ and
$U(1)^r$ for any Kac-Moody algebra  with
$r\times r$ indecomposable generalised Cartan matrix $N$.

In order to motivate the type of representations we will be interested in,
observe that the  homomorphisms $R_n$ above can be
very simply expressed in terms of
 (i) the value of $R$ on the Cartan subalgebra $\langle H\rangle $ (i.e
$R(H)=2z\frac{d}{dz})$
and (ii) a function which transforms as a simple  root vector (i.e $z$
which satisfies $2z\frac{d}{dz}(z)=2z)$.
Let $\hat{\mathfrak{g}}(N)$ be the Kac-Moody algebra (see \ref{theo:Weyl}
below for details) generated by $\{
H_a,X_{\pm i}: 1\le a\le r+s,1\le i\le r\}$  with indecomposable
generalised Cartan matrix $N$  and suppose  given
(i) a homomorphism $F:\langle  H_a: 1\le a\le r+s \rangle\rightarrow
Der({\cal A})$ and (ii) invertible elements
$v_1,\dots,v_r\in {\cal A}$ transforming as simple  root vectors (here
${\cal A}$ is a complex commutative algebra).
The problem we consider in this paper is: can  $F$ be extended to a
homomorphism $\hat{F}:\hat{\mathfrak{g}}(N)\rightarrow
Der({\cal A})$ such that the $\hat{F}(X_{\pm i})$ are of the form $v_i^{\pm
1}\delta_{\pm i}$ with $\delta_{\pm i}\in Im\,F$?
Our main result is the following:
\begin{theobis}
(a) $F$ can be extended as above iff Ker\,$F=Z$ (the centre of
$\hat{\mathfrak{g}}(N)$) and $N$ is the generalised Cartan
matrix of either $sl(r+1,{\bb C})$ or
$\tilde{sl}(r,{\bb C})$ (in  Kac notation $A_r$ or $A_{r-1}^{(1)})$.

(b) In these cases, the extensions of $F$ are in one to one correspondence
with the set of $r\times r$ matrices $A=[A_{ij}]$
satisfying $[\frac{A_{ij}}{A_{jj}}]+[\frac{A_{ji}}{A_{ii}}]=N$
and $\frac{A_{ij}}{A_{jj}}=0$ or $-1$ if $i\not= j$..
\end{theobis}
We show further that there is a ${{\bb Z}^*}^r$-parameter family of
homomorphisms
$\hat{F}^A_{n_1,\dots ,n_r}:\hat{\mathfrak{g}}(N)\rightarrow Der({\cal A})$
associated to each such extension $\hat{F}^A$. The
relationship between
$\hat{F}^A_{n_1,\dots ,n_r}$ and  $\hat{F}^A$ is analogous to the
relationship between $R_n$ and $R$.

 This theorem enables us (cf section 4) to prove that there
exist  families of embeddings of
$sl(r+1,{\bb C})$ and
$\tilde{sl}(r,{\bb C})/Z$ into the Lie algebra of (regular) vector fields
on ${{\bb C}^*}^r$ and  $U(1)^r$ which generalise the
embeddings of $sl(2,{\bb C})$ given above. Explicit formulae in the cases
$sl(2,{\bb C}), \tilde{sl}(2,{\bb
C})/Z,sl(3,{\bb C})$ and $\tilde{sl}(3,{\bb C})/Z$ are given in section 4.

The authors thank D. Bennequin and R.J. Stanton for  useful conversations
and remarks.

\section{Preliminaries}

As general references for this section one can consult the book \cite{Kac}
of V.Kac and the review article \cite{Mac} of I.G.
Macdonald.

Let $N=[n(i,j)]_{1\le i,j\le r}$ be a generalised Cartan matrix
( i.e. $n(i,j)\in \bb Z$, $n(i,i)=2$, $n(i,j)\le 0$ if $i\not= j$,
$n(i,j)=0$ iff $n(j,i)=0$)
and suppose corank$(N)=s$. Consider $3r+s$ independent
variables $\{ H_a,X_{i},X_{-i}: 1\le a\le r+s,1\le i,j\le r\}$, set
\begin{align}
\mathfrak{h}'&=\oplus_{1\le a\le r}\bb C H_a \\
\mathfrak{h}&=\oplus_{1\le a\le r+s}\bb C H_a \\
{\bf  r_+}&=\oplus_{1\le i\le r}\bb C X_i \\
{\bf  r_-}&=\oplus_{1\le i\le r}\bb C X_{-i},
\end{align}
and suppose  $\alpha_1,\dots ,\alpha_r\in\mathfrak{h}^*$ are linearly
independent and satisfy
\begin{equation}\label{alpha}
\alpha_j(H_i)=n(i,j)\qquad \forall 1\le i,j\le r.
\end{equation}
Then $Z=\text{Ker}\alpha_1 \cap\dots \cap\text{Ker}\alpha_r$ is of
dimension $s$, contained in
$\mathfrak{h}' $ and  independent of the choice of $\alpha_1,\dots
,\alpha_r\in\mathfrak{h}^*$ satisfying \eqref{alpha}.

To
the above data one can associate the following Lie algebra ( \cite{Mo} ):
\begin{defi} \label{theo:Weyl}
$\hat{\mathfrak{g}}(N,\alpha_1,\dots ,\alpha_r)$ is the Lie algebra
generated by the $3r+s$ independent
variables $\{ H_a,X_{i},X_{-i}: 1\le a\le r+s,1\le i\le r\}$ subject to the
relations:

\medskip
(a) $[H_a,H_b]=0$;

\medskip
(b) $[X_{i},X_{-i}]=H_i$ and $[X_{i},X_{-j}]=0$ if $i\not=j$;

\medskip
(c) $[H_a,X_j]=\alpha_j(H_a)X_j$ and $[H_a,X_{-j}]=-\alpha_j(H_a)X_{-j}$;

\medskip
(d) $ad^{-n(i,j)+1}(X_i)(X_j)=0$ if $i\not=j$;

\medskip
(e) $ad^{-n(i,j)+1}(X_{-i})(X_{-j})=0$ if $i\not=j$.
\end{defi}
It is well-known ( \cite{Mo} ) that, up to isomorphism, the Lie algebra
$\hat{\mathfrak{g}}(N,\alpha_1,\dots ,\alpha_r)$ is
independent of the choice of $\alpha_1,\dots ,\alpha_r\in\mathfrak{h}^*$
and that its centre is $Z$. Henceforth, we write
$\hat{\mathfrak{g}}(N)$ for $\hat{\mathfrak{g}}(N,\alpha_1,\dots ,\alpha_r)$.

It is also well-known  that
$\hat{\mathfrak{g}}(N)$ is finite-dimensional if and only if all principal
minors of $N$ are $>0$ and then
$\hat{\mathfrak{g}}(N)$ is semisimple ( \cite{Se} ). In this case $s=0$ and
$N$ is  the  Cartan matrix $n(i,j)=\alpha_j(H_i)$
corresponding to a choice of  simple roots $\alpha_1,\dots ,\alpha_r$  and
coroots $H_1,\dots ,H_r$ ( \cite{Se} ).

One can also associate to the above data the  Kac-Moody algebra
$\mathfrak{g}(N)$ ( \cite{Kac} ). This is a quotient of
$\hat{\mathfrak{g}}(N)$ and whether or not the two are isomorphic in
general is an open question. However if $N$ is
symmetrizable then $\hat{\mathfrak{g}}(N)$ is  isomorphic to the Kac-Moody
algebra $\mathfrak{g}(N)$ ( \cite{GaKac} ).

By definition, $\hat{\mathfrak{g}}(N)$ has the following universal property:
\begin{prop} \label{ext}
Let $\mathfrak{l}$ be a Lie algebra. A linear map $F:\mathfrak{h}\oplus
{\bf  r_+}\oplus {\bf  r_-}\rightarrow
\mathfrak{l}$  such that $\{ F(H_a),F(X_i),F(X_{-i}):1\le a\le r+s, 1\le
i\le r \}$ satisfy the relations of \ref{theo:Weyl}
extends to a unique Lie algebra homomorphism
$\hat{F}:\hat{\mathfrak{g}}(N)\rightarrow\mathfrak{l}$.
\end{prop}

\section{Statement of the problem}\label{problems}
Let  ${\cal A}$ be a complex commutative algebra and let
$$
Der({\cal A})=\{D\in End({\cal A}):D(ab)=D(a)b+aD(b)\quad\forall a,b\in
{\cal A}\}.
$$
We suppose given:

\medskip

(A) a Lie algebra homomorphism $F:\mathfrak{h}\rightarrow Der({\cal A})$;

(B) invertible elements $v_{\alpha_1},v_{\alpha_2},\dots ,v_{\alpha_r}$ in
${\cal A}$ such that
$$
F(H)(v_{\alpha_i})=\alpha_i(H)v_{\alpha_i}.
$$
\begin{ex}
Let $L_r={\bb C}[z_1,\frac{1}{z_1},\dots ,z_r,\frac{1}{z_r}]$ be the
algebra of Laurent polynomials in
$z_1,\dots,z_r$. This is the algebra of regular functions on ${{\bb
C}^*}^r$ and $Der(L_r)$ is the Lie algebra of regular
vector fields on  ${{\bb C}^*}^r$.  Define $F:\mathfrak{h}\rightarrow
Der(L_r)$ by setting
$$
F(H)(z_i)=\alpha_i(H)z_i,\quad F(H)(\frac{1}{z_i})=-\alpha_i(H)\frac{1}{z_i},
$$
and extending by the derivation rule (cf \eqref{explicit} for the explicit
formula). Then $F$ and $z_1,\dots ,z_r$ satisfy (A)
and (B).
\end{ex}
\begin{rema}\label{kernel}
By (B), Ker$\,F\subseteq Z=\text{Ker}\alpha_1 \cap\dots
\cap\text{Ker}\alpha_r$ and the $\alpha_i$ are well-defined on
Im$\,F$.
\end{rema}
\begin{rema}
$v_{\alpha_1},\dots ,v_{\alpha_r}$  are linearly independent since
$\alpha_1,\dots
,\alpha_r$ are distinct.
\end{rema}
\begin{rema}
If $\lambda_1,\dots ,\lambda_r\in {\bb C}^*$, then  $F$ and
$\lambda_1v_{\alpha_1},\dots ,\lambda_rv_{\alpha_r}$  satisfy
(A) and (B) above.
\end{rema}

\medskip
In the rest of the article we will usually write $H$ for $F(H)$ to ease
notation.

\medskip
If $\alpha=n_1\alpha_1+n_2\alpha_2+\dots +n_r\alpha_r$ is a $\bb Z$-linear
combination of the $\alpha_i$, we set
$$
v_{\alpha}=v_{\alpha_1}^{n_1}v_{\alpha_2}^{n_2}\dots v_{\alpha_r}^{n_r}.
$$
By the derivation rule it is clear that for all $H\in \mathfrak{h}$
$$
H(v_{\alpha})=\alpha(H)v_{\alpha}.
$$
The commutator of $v_{\alpha}H$ and $v_{\beta}H'$ is given by
\begin{equation}
ad(v_{\alpha}H)(v_{\beta}H')=[v_{\alpha}H,v_{\beta}H']=
v_{\alpha+\beta}(\beta(H)
H'-\alpha(H')H),
\end{equation}
and a  simple induction shows that for $k\ge 2$,
\begin{equation}\label{ad}
\begin{split}
ad^k(v_{\alpha}H)(v_{\beta}H')&=v_{\beta+k\alpha}\bigl(
\beta(H)(\beta+\alpha)(H)\dots(\beta+(k-1)\alpha)(H)\,H'\\
&\qquad -
k\alpha(H')\beta(H)(\beta+\alpha)(H)\dots(\beta+(k-2)\alpha)(H)H\bigr).
\end{split}
\end{equation}

The problem we adress is:

\medskip
{\bf PROBLEM 1}: {\it Is it possible to find  $2r$ elements
$\delta_1,\delta_2\dots\delta_r,\delta_{-1},\delta_{-2}\dots\delta_{-r}$ in
Im$\,F$ such that
$$
\{H_a, X_i=v_{\alpha_i}\delta_{i},X_{-i}=v_{-\alpha_i}\delta_{-i}:1\le a\le
r+s,1\le i\le r\}
$$
satisfy the relations of \ref{theo:Weyl}?}

In fact it will turn out that solving Problem 1 is equivalent to solving
the following a priori simpler problem:

\medskip
{\bf PROBLEM 2}: {\it Is it possible to find  $2r$ elements
$\delta_1,\delta_2\dots\delta_r,\delta_{-1},\delta_{-2}\dots\delta_{-r}$ in
Im$\,F$ such that
$$
\{H_a, X_i=v_{\alpha_i}\delta_{i},X_{-i}=v_{-\alpha_i}\delta_{-i}:1\le a\le
r+s,1\le i\le r\}
$$
satisfy the relations (a),(b) and (c)  of \ref{theo:Weyl}?}

\medskip
The relations (a) and (c) are  satisfied by construction so one only has
consider the relations (b),(d) and (e)
for Problem 1 and (b) for Problem 2. Note that the $\{H_a:r+1\le a\le
r+s\}$ are not involved in these relations.
By \ref{ext}, if there is a solution to Problem 1 then $H_a, X_i,X_{-i}$
generate a Lie subalgebra of $Der({\cal A})$
isomorphic to a quotient of $\hat{\mathfrak{g}}(N)$.

\subsection{Necessary conditions}
We first show that the existence of a solution to Problem 2  imposes
restrictions on the generalised Cartan matrix $N$ and the map
$F:\mathfrak{h}\rightarrow Der({\cal A})$.
\begin{theo}
\label{theo:nec}
Suppose there exist elements
$\delta_1,\delta_2\dots\delta_r,\delta_{-1},\delta_{-2}\dots\delta_{-r}$ in
Im$\,F$ such that
$$
\{H_a, X_i=v_{\alpha_i}\delta_{i},X_{-i}=v_{-\alpha_i}\delta_{-i}:1\le a\le
r+s,1\le i\le r\}
$$
satisfy the relations (b) of \ref{theo:Weyl}. Set
$A_{ij}=\alpha_{i}(\delta_j)$ for $1\le i,j\le r$  (this makes sense by
\ref{kernel}). Then:

\medskip
(i) $A_{ii}\not= 0$ for $1\le i\le r$.

\medskip
(ii) $\delta_{-i}=\frac{1}{A_{ii}}\bigl(-H_i+ \frac{1}{A_{ii}}\delta_{i}
\bigr)$ for $1\le i\le r$.

\medskip
(iii) $\frac{A_{ij}}{A_{jj}}=0$ or $-1$ if $i\not= j$.

\medskip
(iv) $\frac{A_{ij}}{A_{jj}}+\frac{A_{ji}}{A_{ii}}=n(j,i)$.

\medskip
(v) The matrix $A'=[\frac{A_{ij}}{A_{jj}}]$ has the property that if the
$(i,j)$ entry is $-1$ then the $i^{th}$ row and
$j^{th}$column of $A'$ have no other entries equal to $-1$ (their other off
diagonal entries are $0$ by (iii)).

\medskip
(vi) If $\frac{A_{ij}}{A_{jj}}=\frac{A_{ji}}{A_{ii}}=-1$ then the
generalised Cartan
matrix contains the factor $\begin{pmatrix}2&-2\\-2&2\\ \end{pmatrix}$ and
$F(H_i)+F(H_j)=0$.

\end{theo}

Proof.

Rewriting the relations \ref{theo:Weyl}(b) in terms of the $\delta_i$ we
get the following identities in Im$\,F$:
\begin{align}
\label{align: yes}
-\alpha_i(\delta_{i})\delta_{-i}-\alpha_i(\delta_{-i})\delta_{i}&=H_i\\
\label{align: no}
-\alpha_j(\delta_{i})\delta_{-j}-\alpha_i(\delta_{-j})\delta_{i}&=0
\quad\text{ (
if $i\not= j$).}
\end{align}
\begin{rema}\label{invo}
Note that $\lambda_1\delta_1,\dots
,\lambda_r\delta_r,\frac{1}{\lambda_1}\delta_{-1},\dots ,
 \frac{1}{\lambda_r}\delta_{-r}$ ( $\lambda_i\not= 0$) and

\noindent$-\delta_{-1},\dots ,-\delta_{-r},-\delta_{1}\dots ,-\delta_{r}$
are also solutions of this system.
\end{rema}
Applying $\alpha_i$ to the first equation (this is legitimate by
\ref{kernel}) gives
\begin{equation}\label{anymatrix1}
-2\alpha_i(\delta_{i})\alpha_i(\delta_{-i})=\alpha_i(H_i)=2,
\end{equation}
and hence
\begin{equation}\label{anymatrix2}
\alpha_i(\delta_{i})\alpha_i(\delta_{-i})=-1,
\end{equation}
whence  $\alpha_i(\delta_{i})\not= 0$. This proves (i).

 Substituting in (\ref{align: yes}) we get the formula (ii) for
$\delta_{-i}$ in terms of $\delta_{i}$ :
\begin{equation}
\label{equ:minus}
\delta_{-i}=\frac{1}{\alpha_i(\delta_{i})}\bigl(-H_i+
\frac{1}{\alpha_i(\delta_{i})}\delta_{i} \bigr).
\end{equation}
Applying $\alpha_i$ and $\alpha_j$ to \eqref{align: no} gives for all $1\le
i\not= j\le r$
\begin{align}
\alpha_i(\delta_{-j}) \bigl( \alpha_j(\delta_{i})+\alpha_i(\delta_{i})
\bigr)&=0\\
\alpha_j(\delta_{i}) \bigl( \alpha_j (\delta_{-j})+\alpha_i(\delta_{-j})
\bigr)&=0,
\end{align}
and since by \eqref{equ:minus} for all $1\le i, j\le r$
\begin{equation}\label{alpham}
\alpha_j(\delta_{-i})=\frac{1}{\alpha_i(\delta_{i})}\bigl(-n(i,j)+
\frac{\alpha_j(\delta_{i})}{\alpha_i(\delta_{i})}  \bigr)
\end{equation}
these equations are equivalent to for all $1\le i\not= j\le r$
\begin{align}
(-n(j,i)+\frac{A_{ij}}{A_{jj}}) (A_{ji}+A_{ii}) &=0\\
A_{ji}( -1-n(j,i)+\frac{A_{ij}}{A_{jj}}\bigr)&=0.
\end{align}
From this it follows that for $1\le i\not= j\le r$, either
\begin{equation}
A_{ji}=0\quad\text{ and }\quad A_{ij}=A_{jj}n(j,i)
\end{equation}
or
\begin{equation}
A_{ji}=-A_{ii}\quad\text{ and }\quad A_{ij}=A_{jj}(1+n(j,i)).
\end{equation}
Hence $\frac{A_{ji}}{A_{ii}}=0$ or $-1$ if $i\not= j$ and
\begin{equation}
\label{C1}
\frac{A_{ij}}{A_{jj}}+\frac{A_{ji}}{A_{ii}}=n(j,i).
\end{equation}
This proves (iii) and (iv) and the generalised
Cartan matrix is  symmetric with off diagonal entries equal to either
$0,-1$ or $-2$.

Substituting \eqref{C1} in \eqref{alpham} gives for $1\le i,j\le r$
\begin{equation}\label{transpose}
\alpha_j(\delta_{-i})=-\frac{A_{ij}}{A_{jj}A_{ii}}
\end{equation}
and hence \eqref{align: no} becomes
\begin{equation}\label{contra}
A_{ji}\bigl( \delta_{-j}-\frac{1}{A_{jj}A_{ii}}\delta_{i}
\bigr)=0\quad\text{ (if $i\not= j$).}
\end{equation}
If we now apply $\alpha_k$ to this equation we get
\begin{equation}\label{rows}
A_{ji}\bigl( \frac{A_{jk}}{A_{kk}}+\frac{A_{ki}}{A_{ii}}
\bigr)=0\quad\text{ ($i\not= j$ and $k$ arbitrary).}
\end{equation}
We already know that either $\frac{A_{ji}}{A_{ii}}=0$ or $-1$ if $i\not=
j$. Suppose $\frac{A_{ji}}{A_{ii}}=-1$. Then we must
have
\begin{equation}
\frac{A_{jk}}{A_{kk}}+\frac{A_{ki}}{A_{ii}} =0\quad\forall 1\le k\le r.
\end{equation}
But if $k\not= i$ and $k\not=j$, $\frac{A_{jk}}{A_{kk}}\le 0$ and
$\frac{A_{ki}}{A_{ii}} \le 0$;  hence
\begin{equation}
A_{jk}=A_{ki}=0\quad \text{ when $\frac{A_{ji}}{A_{ii}}=-1$, $k\not= i$ and
$k\not=j$}.
\end{equation}
This means that the matrix $A'=[\frac{A_{ij}}{A_{jj}}]$ has the property
that if an entry $\frac{A_{ij}}{A_{jj}}= -1$ then the
row and the column of $A'$ containing that entry have no other entries
equal to $-1$ (i.e., their other off diagonal entries
are $0$). This proves (v).

If $n(i,j)=-2$, we must have $\frac{A_{ij}}{A_{jj}}= -1$ and
$\frac{A_{ji}}{A_{ii}}= -1$ by (iv). By (v) and (iv) the principal
submatrix corresponding to $(i,j)$ defines a decomposition of the
generalised Cartan matrix and the factor corresponding to
$(i,j)$ is the matrix $\begin{pmatrix}2&-2\\-2&2\\ \end{pmatrix}$. This
proves the first part of (vi).

Moreover, by
\eqref{contra}
\begin{equation}
\delta_{-j}-\frac{1}{A_{jj}A_{ii}}\delta_{i}=0\quad\text{ and
}\quad\delta_{-i}-\frac{1}{A_{jj}A_{ii}}\delta_{j}=0,
\end{equation}
which by \eqref{equ:minus} implies
\begin{equation}\label{centre}
H_j=\frac{1}{A_{jj}}\delta_{j}-\frac{1}{A_{ii}}\delta_{i}\quad\text{ and }\quad
H_i=\frac{1}{A_{ii}}\delta_{i}-\frac{1}{A_{jj}}\delta_{j}.
\end{equation}
From this it follows that
\begin{equation}
H_i+H_j=0,
\end{equation}
or more precisely $F(H_i+H_j)=0$ since $H_i$ was an abbreviated notation
for the image of $H_i\in\mathfrak{h}$ under the given
homomorphism $F:\mathfrak{h}\rightarrow Der({\cal A}).$

QED
\begin{cor}\label{cor:nec}
If $N$ is indecomposable and there is a solution to Problem 2 then the Lie
algebra $\hat{\mathfrak g}(N)$ is isomorphic to
either
$A_r$  or $A_{r-1}^{(1)}$.
\end{cor}
Proof. This is a consequence of the following lemma:
\begin{lem}\label{rightcart}
Let $N$ be an $r\times r$ indecomposable generalised Cartan matrix. There
exists an $r\times r$ matrix $A'$ such that

\medskip
(a) $A'_{ii}= 1$ for $1\le i\le r$

\medskip
(b) $A'_{ij}=0$ or $-1$ if $i\not= j$

\medskip
(c) $A'_{ij}+A'_{ji}=n(j,i)$

\medskip
(d) if the $(i,j)$ entry of $A'$ is $-1$ then the $i^{th}$ row and
$j^{th}$column of $A'$ have no other entries equal to $-1$

\medskip
iff $N$ is the generalised Cartan matrix of $A_r\,(r\ge 1)$ or
$A^{(1)}_{r-1}\,(r\ge 2)$. For $A_1$ or $A^{(1)}_1$ there is one
matrix satisfying these conditions and for $A_r\,(r\ge 2)$ or
$A^{(1)}_{r-1}\,(r\ge 3)$ there are two.
\end{lem}
Proof. Suppose there exists a matrix $A'$ with the above properties.  Then
by (b) and (c), $N$ is symmetric and the off
diagonal entries of $N$ are either $0,-1$ or $-2.$

If $n(i,j)=-2$ for some $(i,j)$ then $N$ contains the factor
$\begin{pmatrix}2&-2\\-2&2\\ \end{pmatrix}$
and since $N$ is indecomposable, $N=\begin{pmatrix}2&-2\\-2&2\\
\end{pmatrix}$ which is the generalised Cartan matrix of
$A^{(1)}_1$. In this case $A'=\begin{pmatrix}1&-1\\-1&1\\ \end{pmatrix}$ is
the only possibility.

Suppose that $N$ has no $-2$ entry. Then the Dynkin diagram of $N$ has $r$
vertices and two vertices $i$ and $j$ are connected
iff $n(i,j)=-1.$ No two vertices are connected by more than one line and
since $N$ is indecomposable the diagram is connected.

If $r=1$, $N$ is clearly the Cartan matrix of $A_1$ and there is only one
matrix $A'$ satisfying (a), (b), (c) and (d).
Suppose $r\ge 2$. The matrix $A'$ enables us to orient each line in the
diagram: we put an arrow from $i$ to $j$ if
$A'_{ij}=-1$ and this is a well-defined orientation since by (c), $A'_{ij}$
and $A'_{ji}$ cannot both be equal to $-1$ or $0$.
But by (d), for each vertex
$i$ there is no more than one arrow ending at $i$ and no more than one
arrow beginning at $i$. Hence the Dynkin diagram of $N$
is globally oriented and  must be the Dynkin diagram of either $A_r\,(r\ge
2)$ or
$A^{(1)}_{r-1}\,(r\ge 3)$. There are then exactly two  matrices satisfying
(a), (b), (c) and (d): $A'$ and
${}^tA'$ (corresponding to the two possible global orientations of the
Dynkin diagram).

For the converse, note that the (generalised) Cartan matrix of $A_r$ can be
written
\begin{equation}\label{Afin}
\begin{pmatrix}2&-1&0&\dots&0\\
-1&2&\ddots&\ddots&\vdots\\
 0&\ddots&\ddots&\ddots&0\\
\vdots&\ddots&\ddots&\ddots&-1\\
0&\dots&0&-1&2\\ \end{pmatrix}=
\begin{pmatrix}1&-1&0&\dots&0\\
0&1&\ddots&\ddots&\vdots\\
0&\ddots&\ddots&\ddots&0\\
\vdots&\ddots&\ddots&\ddots&-1\\
0&\dots&0&0&1\\ \end{pmatrix}
+
\begin{pmatrix}1&0&0&\dots&0\\
-1&1&\ddots&\ddots&\vdots\\
0&\ddots&\ddots&\ddots&0\\
\vdots&\ddots&\ddots&\ddots&0\\
0&\dots&0&-1&1\\ \end{pmatrix}
\end{equation}
and that the generalised Cartan matrix of $A_{r-1}^{(1)}, r\ge 2$ can be
written

\begin{equation}\label{Ainfin}
\begin{pmatrix}2&-1&0&\dots&-1\\
-1&2&\ddots&\ddots&0\\
 0&\ddots&\ddots&\ddots&\vdots\\
\vdots&\ddots&\ddots&\ddots&-1\\
-1&\dots&0&-1&2\\ \end{pmatrix}=
\begin{pmatrix}1&-1&0&\dots&0\\
0&1&\ddots&\ddots&\vdots\\
0&\ddots&\ddots&\ddots&0\\
\vdots&\ddots&\ddots&\ddots&-1\\
-1&\dots&0&0&1\\ \end{pmatrix}
+
\begin{pmatrix}1&0&0&\dots&-1\\
-1&1&\ddots&\ddots&\vdots\\
0&\ddots&\ddots&\ddots&0\\
\vdots&\ddots&\ddots&\ddots&0\\
0&\dots&0&-1&1\\ \end{pmatrix}
\end{equation}
The matrices on the righthand sides of these equations clearly satisfy (a),
(b), (c) and (d).

QED

For future reference we introduce the following terminology:
\begin{defi}\label{smdef}

(i) Let $N$ be the generalised Cartan matrix of $A_r$ or $A^{(1)}_{r-1}$.
An $r\times r$ matrix $A=[A_{ij}]$ is called a
solution matrix of $N$ if the matrix $[\frac{A_{ij}}{A_{jj}}]$ satisfies
(a), (b), (c) and (d) of \ref{rightcart}. If in
addition  $A_{ii}=1$ for all $1\le i\le r$, $A$ is called a normalised
solution matrix of $N$.

(ii) If $A$ is a solution matrix  we set
$I^A=\{(i,j):\frac{A_{ji}}{A_{ii}}=-1\}$. Note that $I^{({}^tA)}={}^t(I^A)$.
\end{defi}
 From  \eqref{Afin} and \eqref{Ainfin} we see that if $A$ is a solution
matrix of
$A_r$ then (up to transposition)
\begin{equation}
I^A=\{(1,2),(2,3),\dots (r-1,r)\}
\end{equation}
and that if $A$ is a solution  matrix of $A^{(1)}_{r-1}$
\begin{equation}\label{incid}
I^A=\{(1,2),(2,3),\dots (r-1,r),(r,1)\}.
\end{equation}
This has the following consequence:
\begin{cor}
If  $N$ is indecomposable and there is a solution to Problem 2 then
Ker$\,F= Z$.
\end{cor}
Proof. By \ref{kernel}, Ker$\,F\subseteq Z$. If $N$ is the generalised
Cartan matrix of $A_r$ then $Z=\{0\}$ and there is
nothing to prove. If
$N$ is the generalised Cartan matrix of $A^{(1)}_{r-1}$ then it is known (
\cite{Kac}, \cite{Mac} ) that $Z\subset\mathfrak{h}$
is spanned by
$H_1+H_2+\dots +H_r.$ But by
\eqref{centre}
\begin{equation}
\forall (i,j)\in I^A,\quad
H_j=\frac{1}{A_{jj}}\delta_{j}-\frac{1}{A_{ii}}\delta_{i}
\end{equation}
and hence by \eqref{incid}
\begin{equation}
H_1+H_2+\dots +H_r=0,
\end{equation}
which, recalling our abbreviated notation, really means $F(H_1+H_2+\dots
+H_r)=0$.
QED
\begin{rema}
The only property of the generalised Cartan matrix $N$ used in this section
was that its diagonal elements are invertible
complex numbers. This was needed (cf \eqref{anymatrix1},\eqref{anymatrix2})
to prove that $A_{ii}\not= 0$.
\end{rema}

\subsection{Equivalence of Problem 1 and Problem 2}
\begin{theo}
\label{theo:rel(d)}
Suppose there exist elements
$\delta_1,\delta_2\dots\delta_r,\delta_{-1},\delta_{-2}\dots\delta_{-r}$ of
Im$\,F$ such that
$$
\{H_a, X_i=v_{\alpha_i}\delta_{i},X_{-i}=v_{-\alpha_i}\delta_{-i}:1\le a\le
r+s,1\le i\le r\}
$$
satisfy the relations (b) of \ref{theo:Weyl}. Then
$$
\{H_a, X_i=v_{\alpha_i}\delta_{i},X_{-i}=v_{-\alpha_i}\delta_{-i}:1\le a\le
r+s,1\le i\le r\}
$$
also satisfy the relations (d) and (e) of \ref{theo:Weyl}.
\end{theo}
\begin{cor}\label{equiv}
Problem 1 has a solution iff Problem 2 has a solution.
\end{cor}
Proof.
To prove that the relations (d) of \ref{theo:Weyl} are satisfied, we have
to show that
$$
ad^{-n(i,j)+1}(X_i)(X_j)=0 \qquad\text{ if $i\not=j$},
$$
which is equivalent to
$$
ad^{-n(i,j)+1}(v_{\alpha_i}\delta_{i})(v_{\alpha_j}\delta_{j})=0\qquad
\text{ if $i\not=j$}.
$$
By \ref{theo:nec}(iii) and (iv), if the relations (b) of \ref{theo:Weyl}
are satisfied then the Cartan matrix is symmetric and
$n(i,j)=0,-1$ or $-2$ if $i\not= j$.

Setting $\alpha=\alpha_i,H=\delta_{i},\beta=\alpha_j$ and $H'=\delta_{j}$
in the formula \eqref{ad} we get
\begin{equation}
\label{equ:ad1}
ad^{-n(i,j)+1}(v_{\alpha_i}\delta_{i})(v_{\alpha_j}\delta_{j})=
v_{\alpha_i+\alpha_j}\bigl(
A_{ji}\delta_{j}-A_{ij}\delta_{i}\bigr)
\end{equation}
if $n(i,j)=0$,
\begin{equation}
\label{equ:ad2}
ad^{-n(i,j)+1}(v_{\alpha_i}\delta_{i})(v_{\alpha_j}\delta_{j})=
v_{\alpha_j+2\alpha_i}\bigl(
A_{ji}(A_{ji}+A_{ii})\delta_{j}-2A_{ij}A_{ji}\delta_{i}\bigr)
\end{equation}
if $n(i,j)=-1$ and
\begin{equation}\label{equ:ad3}
\begin{split}
ad^{-n(i,j)+1}(v_{\alpha_i}\delta_{i})(v_{\alpha_j}\delta_{j})&=
v_{\alpha_j+3\alpha_i}\bigl(
A_{ji} (A_{ji}+A_{ii}) (A_{ji}+2A_{ii})\delta_{j}\\
&\qquad -
3A_{ij}A_{ji}(A_{ji}+A_{ii})\delta_{i}\bigr).
\end{split}
\end{equation}
if $n(i,j)=-2$.

In the first case we  have $A_{ij}=A_{ji}=0$ by \ref{theo:nec}(iii) and
\ref{theo:nec}(iv) and hence by \eqref{equ:ad1}
$$
ad^{-n(i,j)+1}(v_{\alpha_i}\delta_{i})(v_{\alpha_j}\delta_{j})=0.
$$

In the second case either $A_{ij}=-A_{jj}$ and $A_{ji}=0$ or $A_{ij}=0$ and
$A_{ji}=-A_{ii}$, again by
\ref{theo:nec}(iii) and \ref{theo:nec}(iv). Thus by \eqref{equ:ad2}
$$
ad^{-n(i,j)+1}(v_{\alpha_i}\delta_{i})(v_{\alpha_j}\delta_{j})=0.
$$
if $n(i,j)=-1$.

In the third case  $A_{ij}=-A_{jj}$ and $A_{ji}=-A_{ii}$  by
\ref{theo:nec}(iii) and \ref{theo:nec}(iv).
Thus by \eqref{equ:ad3},
$$
ad^{-n(i,j)+1}(v_{\alpha_i}\delta_{i})(v_{\alpha_j}\delta_{j})=0
$$
if $n(i,j)=-2$. Hence in all cases the relations (d) of \ref{theo:Weyl} are
satisfied.

To prove that the relations (e) of \ref{theo:Weyl} are satisfied we have to
show that
$$
ad^{-n(i,j)+1}(v_{-\alpha_i}\delta_{-i})(v_{-\alpha_j}\delta_{-j})=0\qquad
\text{ if $i\not=j$}.
$$
Setting $B_{ij}=-\alpha_i(\delta_{-j})$, by \eqref{ad} this is equivalent to
\begin{equation}\label{Bad1}
B_{ji}\delta_{-j}-B_{ij}\delta_{-i}=0
\end{equation}
if $n(i,j)=0$, to
\begin{equation}\label{Bad2}
B_{ji}(B_{ji}+B_{ii})\delta_{-j}-2B_{ij}B_{ji}\delta_{-i}=0
\end{equation}
if $n(i,j)=-1$ and to
\begin{equation}\label{Bad3}
B_{ji} (B_{ji}+B_{ii}) (B_{ji}+2B_{ii})\delta_{-j}-
3B_{ij}B_{ji}(B_{ji}+B_{ii})\delta_{-i}=0
\end{equation}
if $n(i,j)=-2$.

But by \eqref{transpose}, we have  forall $1\le i,j\le r$,
\begin{equation}
\frac{B_{ji}}{B_{ii}}=\frac{A_{ij}}{A_{jj}}.
\end{equation}
This means the $B_{ij}$ satisfy the properties (i), (iii) and (iv) of
\ref{theo:nec} and the equations
\eqref{Bad1},\eqref{Bad2} and \eqref{Bad3} follow from them exactly as the
equations
\eqref{equ:ad1},\eqref{equ:ad2} and \eqref{equ:ad3} followed from the
properties of the $A_{ij}$ in the above proof of the
relations (d) of \ref{theo:Weyl}.

QED

\subsection{Sufficient conditions}\label{ansatz}

Let $N$ be the Cartan matrix of  either $A_r\,(r\ge 1)$  or
$A_{r-1}^{(1)}\,( r\ge 2)$ and let $A=[A_{ij}]$ be a solution
matrix of $N$ (cf \ref{smdef}). We suppose given a complex commutative
algebra ${\cal A}$ and:

(A) a Lie algebra homomorphism $F:\mathfrak{h}\rightarrow Der({\cal A})$
such that Ker\,$F=Z$;

(B) invertible elements $v_{\alpha_1},v_{\alpha_2},\dots ,v_{\alpha_r}$ in
${\cal A}$ such that
$$
F(H)(v_{\alpha_i})=\alpha_i(H)v_{\alpha_i}.
$$
As before,  we will write $H$ for $F(H)$ to ease notation.

\medskip
By hypothesis $F$ factors to an isomorphism $F:\mathfrak{h}\,/Z\rightarrow$
Im$\,F$ and therefore $\{\alpha_1,\dots
,\alpha_r\}$ factors to a basis of (Im$\,F)^*$. Hence there exist unique
$\delta_1,\delta_2\dots\delta_r$  in Im$\,F$ such
that
$[\alpha_i(\delta_j)]=[A_{ij}]$. For $1\le i\le r$ set
\begin{equation}\label{here}
\delta_{-i}=\frac{1}{A_{ii}}\bigl(-H_i+  \frac{1}{A_{ii}} \delta_{i} \bigr).
\end{equation}
\begin{theo}
$\{ H_i, X_i=v_{\alpha_i}\delta_{i},X_{-i}=v_{-\alpha_i}\delta_{-i}:1\le
i\le r \}$ satisfy the relations (b) of
\ref{theo:Weyl} and hence Problem 1 has a solution.
\end{theo}

Proof. By \eqref{align: yes} and \eqref{align: no}, the relations (b) of
\ref{theo:Weyl} are satisfied iff
\begin{align}
\label{align: yess}
-\alpha_i(\delta_{i})\delta_{-i}-\alpha_i(\delta_{-i})\delta_{i}&=H_i\\
\label{align: noo}
-\alpha_j(\delta_{i})\delta_{-j}-\alpha_i(\delta_{-j})\delta_{i}&=0
\quad\text{ (
if $i\not= j$).}
\end{align}
Equation \eqref{align: yess} is an immediate consequence of \eqref{here}.
Also by \eqref{here},  \eqref{align: noo} is
equivalent to
\begin{equation}
A_{ji}\bigl( H_j-\frac{\delta_{j}}{A_{jj}}+\frac{\delta_{i}}{A_{ii}}
\bigr)=0\quad\text{ (if $i\not= j$)},
\end{equation}
which, since  $(i,j) \in I^A\Leftrightarrow \frac{A_{ji}}{A_{ii}}=-1$ and
otherwise $A_{ji}=0$, is
equivalent to
\begin{equation}
\forall (i,j) \in I^A,\qquad \qquad
H_j-\frac{\delta_{j}}{A_{jj}}+\frac{\delta_{i}}{A_{ii}}=0.
\end{equation}
The lefthand side  is in Im$\,F$ so this equation is equivalent to
\begin{equation}
\forall (i,j) \in I^A,\forall k: 1\le k\le
r,\quad\alpha_k(H_j-\frac{\delta_{j}}{A_{jj}}+\frac{\delta_{i}}{A_{ii}})=0
\end{equation}
since $\alpha_1,\dots ,\alpha_r$ is a basis of (Im$\,F)^*.$
But for all $1\le i,j,k\le r$ we have
\begin{align}
\alpha_k(H_j-\frac{\delta_{j}}{A_{jj}}+\frac{\delta_{i}}{A_{ii}})&=n(j,k)-
\frac{A_{kj}}{A_{jj}}+\frac{A_{ki}}{A_{ii}}\\
&=\frac{A_{jk}}{A_{kk}}+\frac{A_{ki}}{A_{ii}}\quad\text{(by
\ref{rightcart}(c))}.
\end{align}
Suppose $\frac{A_{ji}}{A_{ii}}=-1$. Then for $k=i$ or $k=j$ this is
obviously $0$  and for $k\not=i,k\not=
j$, both terms vanish by property \ref{rightcart}(d). This proves the theorem.

QED

In conclusion:
\begin{cor}
If $N$ is the generalised Cartan matrix of $A_r$ or $A^{(1)}_{r-1}$ and
$F:\mathfrak{h}\rightarrow Der({\cal A})$ satisfies
Ker$\,F=Z$ then Problem 1 has a solution. The set of solutions is in one to
one correspondence with the set of solution
matrices of $N.$
\end{cor}

\begin{rema}
 If $[A_{ij}]$ is a solution matrix (SM), $[\mu_jA_{ij}]$ is also a SM
where the $\mu_j$ are arbitrary
non-zero complex numbers. This gives a principal ${\bb C^*}^r$ action on
the space of SMs and the quotient is the space of
normalised SMs. One clearly has an isomorphism of principal ${\bb C^*}^r$
bundles:
$$
\text{Solution Matrices}\cong \text{Normalised Solution Matrices}\times
{\bb C^*}^r\
$$
given by:
$$
[A_{ij}]\rightarrow \bigl( [\frac{A_{ij}}{A_{jj}}],(A_{11},A_{22},\dots
,A_{rr})\bigr).
$$
In the case of $A_1$ and $A_{1}^{(1)}$ there is one  normalised solution
matrix which is symmetric. In the other cases there
are two,  one being the transpose of the other (cf \ref{rightcart}).

The  space of SMs has a natural involution $[A_{ij}]\mapsto [A_{ji}]$ which
commutes with the principal ${\bb C^*}^r$ action.
In terms of solutions $(\delta_1,\dots\delta_r,\delta_{-1},\dots
,\delta_{-r})$ this corresponds to (cf. \ref{invo})
$$
(\delta_1,\dots\delta_r,\delta_{-1},\dots ,\delta_{-r})\mapsto
(-\delta_{-1},\dots ,-\delta_{-r},-\delta_{1}\dots ,-\delta_{r}).
$$
\end{rema}
\begin{rema} \label{discreteparam}
Since $Z=\cap_i\,Ker \alpha_i=\langle H_1+\dots +H_r\rangle$, we can find a
basis  $H_1',\dots,H_r',H_1+\dots +H_r$ of
$\mathfrak{h}$ such that
$\alpha_i(H_j')=\delta_{ij}$ for
$1\le i,j\le r$. For $(n_1,\dots ,n_r)\in {{\bb
Z}^*}^r$,  define $F_{n_1,\dots ,n_r}:\mathfrak{h}\rightarrow Der({\cal A})$ by
\begin{align}
F_{n_1,\dots ,n_r}(H_i')&=\frac{1}{n_i}F(H_i')\\
F_{n_1,\dots ,n_r}(H_1+\dots +H_r)&=0.
\end{align}
One checks  that $F_{n_1,\dots ,n_r}$ is independent of the choice of
basis and that $F_{n_1,\dots ,n_r}$ and
$v_1^{n_1},\dots ,v_r^{n_r}\in {\cal A}$ satisfy properties (A) and (B)
above. Hence each  SM defines an extension
of $F_{n_1,\dots ,n_r}$ to a homomorphism $\hat{F}^A_{n_1,\dots
,n_r}:\hat{\mathfrak{g}}(N)\rightarrow Der({\cal A})$.
\end{rema}

\section{Examples: $sl(r),\tilde{sl}(r)$ and Lie algebras of vector fields
on tori }

 In this section $N$ is the  generalised Cartan matrix of  either
$A_r\,(r\ge 1)$ or $A^{(1)}_{r-1}\,(r\ge 2)$. The Lie
algebra $\hat{\mathfrak{g}}(N)$ is therefore isomorphic to either $A_r$ or
$A^{(1)}_{r-1}$ (i.e to $sl(r+1)$ or
$\tilde{sl}(r)$).

If $M$ is a manifold, $C^{\infty}(M)$,the space of smooth complex-valued
functions on $M$, is a commutative algebra and
$Der\bigl(C^{\infty}(M)\bigr)$ is isomorphic to the Lie algebra of vector
fields ${\cal X}(M)$ on $M$. Suppose that on $M$ one
can find invertible functions $v_1,\dots ,v_r$ and  commuting vector fields
$D_1,\dots ,D_r$ such that
\begin{equation}\label{orthogonal}
D_i(v_j)=\delta_{ij}v_j.
\end{equation}
Then if  we define $F:\mathfrak{h}\rightarrow {\cal X}(M)$ by
\begin{equation}\label{explicit}
F(H_a)=\sum_{k=1}^{k=r}\alpha_k(H_a)D_k\quad (1\le a\le r+s),
\end{equation}
it is clear that Ker$\,F=Z$ and that
$$
F(H_a)(v_j)=\alpha_j(H_a)v_j.
$$
Hence $F$  satisfies the hypotheses (A) and (B) of subsection \ref{ansatz} and
  each choice of a solution matrix $A$ defines
an extension of this map to a Lie algebra homomorphism
$\hat{F}^A:\hat{\mathfrak{g}}(N)\rightarrow{\cal X}(M)$.

$A_r\,(r\ge 1)$ is simple so $\hat{F}^A$ is injective in that case. The Lie
algebra  $A^{(1)}_{r-1}\,(r\ge 2)$ is not simple
but it is known (Proposition 1.7(b) in \cite{Kac}) that every ideal either
contains the derived algebra ${A^{(1)}_{r-1}}'$ or
is contained in the centre $Z$. It is clear that Ker$\,\hat{F}^A$ does not
contain ${A^{(1)}_{r-1}}'$ and therefore
Ker$\,\hat{F}^A=Z$. Hence
$\hat{F}^A:\hat{\mathfrak{g}}(N)\rightarrow{\cal X}(M)$ factors to an
injection of $A^{(1)}_{r-1}/Z$ into ${\cal X}(M)$.

To describe $\hat{F}^A$ more explicitly one has to calculate the $\delta_i$
and $\delta_{-i}$ associated to a solution matrix
$[A_{ij}]$ (cf. subsection \ref{ansatz}). From equation \eqref{explicit} it
is easy to see that
$$
\alpha_i(D_j)=\delta_{ij}\quad\forall 1\le i,j\le r,
$$
and hence that for $1\le i\le r$,
\begin{equation}
\delta_i=\sum_{j=1}^{j=r}A_{ji}D_j \quad\text{and}\quad
\delta_{-i}=-\sum_{j=1}^{j=r}\frac{A_{ij}}{A_{ii}A_{jj}}D_j.
\end{equation}
By subsection \ref{ansatz} this proves the
\begin{prop}\label{expemb}
With the notation above, the $\bb C$- linear map ${F}^A:\mathfrak{h}\oplus
{\bf  r_+}\oplus {\bf  r_-}\rightarrow
{\cal X}(M)$ given by
\begin{align}
{F}^A(H_a)&=\sum_{j=1}^{j=r}\alpha_j(H_a)D_j\quad ( 1\le a\le r+s)\\
{F}^A(X_i)&=A_{ii}v_i\sum_{j=1}^{j=r}\frac{A_{ji}}{A_{ii}}D_j\quad ( 1\le
i\le r)\\
{F}^A(X_{-i})&=-\frac{1}{A_{ii}v_i}\sum_{j=1}^{j=r}\frac{A_{ij}}{A_{jj}}D_j
\quad( 1\le i\le r)
\end{align}
extends to a Lie algebra homomorphism
$\hat{F}^A:\hat{\mathfrak{g}}(N)\rightarrow{\cal X}(M)$ with kernel $Z$.
\end{prop}
If $n_1,\dots ,n_r\in{\bb Z}^*$, the vector fields $\frac{D_1}{n_1},\dots
\frac{D_r}{n_r}$ and the functions
$v_1^{n_1},\dots ,v_r^{n_r}$ also satisfy \eqref{orthogonal} (cf
\ref{discreteparam}) and hence, repeating the above
calculations in this case, we get the
\begin{cor}\label{expembn}
Let $n_1,\dots ,n_r$ be non-zero integers.The $\bb C$- linear map
${F}^A_{n_1,\dots ,n_r}:\mathfrak{h}\oplus
{\bf  r_+}\oplus
{\bf  r_-}\rightarrow
{\cal X}(M)$ given by
\begin{align}
{F}^A_{n_1,\dots
,n_r}(H_a)&=\sum_{j=1}^{j=r}\alpha_j(H_a)\frac{D_j}{n_j}\quad ( 1\le a\le
r+s)\\
{F}^A_{n_1,\dots
,n_r}(X_i)&=A_{ii}v_i^{n_i}\sum_{j=1}^{j=r}\frac{A_{ji}}{A_{ii}}
\frac{D_j}{n_j}\quad ( 1\le i\le r)\\
{F}^A_{n_1,\dots
,n_r}(X_{-i})&=-\frac{1}{A_{ii}v_i^{n_i}}\sum_{j=1}^{j=r}\frac{A_{ij}}
{A_{jj}}\frac{D_j}{n_j}
\quad ( 1\le i\le r)
\end{align}
extends to a Lie algebra homomorphism $\hat{F}^A_{n_1,\dots,n_r}:
\hat{\mathfrak{g}}(N)\rightarrow{\cal X}(M)$ with kernel $Z$.
\end{cor}

In particular, taking $M={\bb C^*}^r$, $v_j=z_j$ (the coordinate functions)
and $D_j=z_j\frac{\partial}{\partial z_j}$ for
$1\le j\le r$, we can construct  homomorphisms from $sl(r+1)$ or
$\tilde{sl}(r)$ to the Lie algebra of regular vector fields on ${\bb
C^*}^r$. Similarly, taking
 $M={U(1)}^r$, $v_j=e^{i\theta_j}$  and $D_j=-i\frac{\partial}{\partial
\theta_j}$ for
$1\le j\le r$, we can construct  homomorphisms from $sl(r+1)$ or
$\tilde{sl}(r)$ to the Lie algebra complex vector fields on the compact torus
${U(1)}^r$. We give now the complete formulae for $\hat{F}^A_{n_1,\dots
,n_r}$ for $sl(2),sl(3),\tilde{sl}(2)$ and
$\tilde{sl}(3)$.
\begin{ex}
If $\mathfrak{g}=sl(2,\bb C)$ and $H_1,X_1,X_{-1}$ is the standard basis, a
solution matrix
$[A_{ij}]$ is  a complex number $\lambda\not= 0$. The embeddings
$\hat{F}^{\lambda}_n:sl(2,\bb C)\rightarrow {\cal X}(M)$
are then given by
\begin{equation}
H_1\mapsto 2\frac{D_1}{n},\quad X_1\mapsto \lambda
v_1^{n}\frac{D_1}{n},\quad X_{-1}\mapsto -\frac{1}{\lambda}
v_1^{-n}\frac{D_1}{n}.
\end{equation}
If $M={\bb C}^*$, $v_1=z$ and $D_1=z\frac{d}{dz}$ this gives the well-known
embeddings:
\begin{equation}
H_1\mapsto \frac{1}{n}2z\frac{d}{dz},\quad X_1\mapsto \frac{\lambda}{n}
z^{n+1}\frac{d}{dz},\quad X_{-1}\mapsto
-\frac{1}{n\lambda} z^{-n+1}\frac{d}{dz}
\end{equation}
of $sl(2,\bb C)$ in ${\cal X}({\bb C^*})$.
The representations $\hat{F}^{\lambda}_{\pm 1}$ can be  obtained
geometrically by restricting the derivative of the action of
$PSL(2,{\bb C})$ on $P_1({\bb C})$ to an appropriate affine coordinate chart.

\end{ex}
\begin{ex} (See also \cite{R}).
If $\mathfrak{g}=sl(3,\bb C)$, let $\{H_1,H_2,X_{\pm i}:1\le i\le 3\}$ be a
standard basis, i.e.: $H_1,H_2$ span a Cartan
subalgebra and are the coroots corresponding to the  simple roots
$\alpha_1,\alpha_2$; the positive roots are
$\alpha_1,\alpha_2,\alpha_3=\alpha_1+\alpha_2$; the Cartan matrix
$N=[\alpha_j(H_i)]_{1\le i,j\le
2}$ is the matrix $\begin{pmatrix} 2&-1\\-1&2\\ \end{pmatrix}$;
$\{H_i,X_{\pm i}:1\le i\le 2\}$ satisfy the relations
of \ref{theo:Weyl} and $X_{\pm 3}=\pm [X_{\pm 1},X_{\pm 2}]$. The solution
matrices of $N$ are of the form
\begin{equation}
A=\begin{pmatrix}
A_{11}&\frac{1}{2}(-1-\varepsilon)A_{22}\\\frac{1}{2}(-1+\varepsilon)
A_{11}&A_{22}\\ \end{pmatrix},
\end{equation}
where $ A_{11}, A_{22}\in \bb C^*$ are arbitrary and $\varepsilon=\pm 1$.
The embeddings $\hat{F}^A_{n_1,n_2}:sl(3,\bb
C)\rightarrow
{\cal X}(M)$  are then given by
\begin{eqnarray}
\label{eq:sl3}
&&\left\{ \begin{array}{l}
X_1 \mapsto A_{11} v_1^{n_1}
\left(\frac{D_1}{n_1} +
 \frac{1}{2}(-1+\varepsilon) \frac{D_2}{n_2}\right)  \cr
X_2 \mapsto A_{22} v_2^{n_2}
\left( \frac{1}{2}(-1-\varepsilon) \frac{D_1}{n_1}+
\frac{D_2}{n_2}\right) \cr
X_3\mapsto A_{11}A_{22} v_1^{n_1 }v_2^{n_2}
\left(-\frac{1}{2}(-1-\varepsilon) \frac{D_1}{n_1}
+\frac{1}{2}(-1+\varepsilon)\frac{D_2}{n_2}\right)
\end{array}
\right.  \nonumber \\
 \nonumber \\
&&\left\{
\begin{array}{l}
X_{-1} \mapsto  \frac{-1}{A_{11}v_1^{n_1}}
\left( \frac{D_1}{n_1}
+ \frac{1}{2}(-1-\varepsilon) \frac{D_2}{n_2}\right)
 \cr
X_{-2}\mapsto  \frac{-1}{A_{22}v_2^{n_2}}
\left( \frac{1}{2}(-1+\varepsilon)  \frac{D_1}{n_1}
 + \frac{D_2}{n_2}\right)\cr
X_{-3}\mapsto \frac{1}{A_{11} A_{22}v_1^{n_1}v_2^{n_2}}
\left(-\frac{1}{2}(-1+\varepsilon) \frac{D_1}{n_1}+
\frac{1}{2}(-1-\varepsilon)
\frac{D_2}{n_2} \right)
\end{array}
\right.
\nonumber \\
  \\
&&\left\{
\begin{array}{l}
H_1 \mapsto  2\frac{D_1}{n_1}-\frac{D_2}{n_2}  \cr
H_2 \mapsto  -\frac{D_1}{n_1}+2 \frac{D_2}{n_2}
\end{array}
\right.
\nonumber
\end{eqnarray}
If $M={{\bb C}^*}^2$, the representations $\hat{F}^{A}_{\pm1,\pm1}$ can be
obtained geometrically by restricting the
derivative of the action of $PSL(3,{\bb C})$ on $P_2({\bb C})$ to an
appropriate affine coordinate chart.

\end{ex}

To give the explicit formulae for the embeddings of the Lie algebras
$\tilde{sl}(2,{\bb C})/Z$ and $\tilde{sl}(3,{\bb C})/Z$
(or $A^{(1)}_1$ and $A^{(1)}_2$ in Kac notation) it is convenient to first
describe them using loop
algebras ( \cite{Kac}, \cite{Mac} ). Recall that in terms of
\ref{theo:Weyl}, $\tilde{sl}(r,{\bb C})$ is constructed from
generators $\{
\tilde{H}_0,\dots ,
\tilde{H}_r,\tilde{X}_{\pm 0},\dots ,\tilde{X}_{\pm (r-1)}\}$, the
generalised Cartan matrix \eqref{Ainfin} and linear forms
$\tilde{\alpha}_0,\dots ,\tilde{\alpha}_{r-1}$ on $\langle
\tilde{H}_0,\dots,\tilde{H}_r\rangle$  whose values on $\langle
\tilde{H}_0,\dots,\tilde{H}_{r-1}\rangle$ are given by the Cartan matrix and
whose values on $\tilde{H}_r$ are given by
$\tilde{\alpha}_i(\tilde{H}_r)=\delta_{i0}$.

The underlying vector space of the loop algebra $L\left(sl(r,{\bb
C})\right)$  is the tensor product of $sl(r,{\bb C})$ with
Laurent series in  $t$ and its bracket is
\begin{equation}
[t^mX,t^nY]=t^{m+n}[X,Y]\quad(m,n\in{\bb Z},X,Y\in sl(r,{\bb C})).
\end{equation}
The operator $d:L\left(sl(r,{\bb C})\right)\rightarrow L\left(sl(r,{\bb
C})\right)$ given by
\begin{equation}
d(t^mX)=mt^{m}X
\end{equation}
is a derivation and one can form the semi-direct product
$\tilde{L}\left(sl(r,{\bb C})\right)=L\left(sl(r,{\bb
C})\right)\oplus {\bb C}d$ with the bracket
\begin{equation}
[\xi\oplus \lambda d,\eta\oplus \mu d]=[\xi,\eta]+\lambda d(\eta)-\mu d(\xi).
\end{equation}
It is well-known ( \cite{Kac}, \cite{Mac} ) that the Lie algebras
$\tilde{sl}(r,{\bb C})/Z$ and $\tilde{L}\left(sl(r,{\bb
C})\right)$ are isomorphic. Explicitly, the map $\Psi:\mathfrak{h}\oplus
{\bf  r_+}\oplus {\bf  r_-}\rightarrow
\tilde{L}\left(sl(r,{\bb C})\right)$ extends to an isomorphism
$\tilde{\Psi}:\tilde{sl}(r,{\bb C})/Z\rightarrow
\tilde{L}\left(sl(r,{\bb C})\right)$ where
\begin{align}\label{loopiso}
\Psi(\tilde{H}_r)&=d\nonumber\\
\Psi(\tilde{H}_{i})&=H_i\quad(\text{ if }1\le i\le r-1)\nonumber\\
\Psi(\tilde{H}_{0})&=-H_{\phi}=-(H_1+\dots +H_{r-1})\nonumber\\
\Psi(\tilde{X}_{\pm i})&=X_{\pm i}\quad(\text{ if }1\le i\le r-1)\nonumber\\
\Psi(\tilde{X}_{\pm 0})&=t^{\pm 1}X_{\mp\phi }.
\end{align}
Here $\{H_1,\dots ,H_{r-1},X_{\pm 1},\dots ,X_{\pm (r-1)}\}$ is a standard
generating set of $sl(r,{\bb C})$ corresponding
to simple  roots $\alpha_1,\dots ,\alpha_{r-1}$,
$\phi=\alpha_1+\dots +\alpha_{r-1}$ is the highest root and $X_{\pm\phi }$
are root vectors satisfying
$[X_{\phi},X_{-\phi }]=H_{\phi}$, the coroot corresponding to $\phi$.
\begin{ex}
If $\{H_1,X_{\pm 1}\}$ is a standard basis of $sl(2,{\bb C})$ then
$\{d,t^mH_1,t^nX_{\pm 1}:m,n\in {\bb Z}\}$
is a basis of $\tilde{L}\left(sl(2,{\bb C})\right)$. The Cartan matrix of
$\tilde{sl}(2,{\bb C})$ is
$\begin{pmatrix} 2&-2\\ -2&2\\\end{pmatrix}$ and the solution matrices are
of the form $\begin{pmatrix} A_{00}&-A_{11}\\
-A_{00}&A_{11}\\\end{pmatrix}$. The embedding $\hat{F}^A_{n_0,n_1}\circ
\Psi^{-1}:\tilde{L}\left(sl(2,{\bb
C})\right)\rightarrow {\cal X}(M)$ is  given by the formulae:
\begin{eqnarray}
\label{eqtildesl2}
&&\begin{array}{l}
t^mX_1 \mapsto    A_{11}\left(A_{00} A_{11} \right)^m
\left(v_0^{n_0} v_1^{n_1} \right)^m {v_1}^{n_1} (-\frac{D_0}{n_0}+
\frac{D_1}{n_1}) \cr
t^mX_{-1} \mapsto - \frac{1}{A_{11}}\left({A_{00} A_{11}} \right)^m
\left(v_0^{n_0} v_1^{n_1} \right)^mv_1^{-n_1} (-\frac{D_0}{n_0}+
\frac{D_1}{n_1}) \cr
t^mH_1 \mapsto   \left(A_{00} A_{11} \right)^m \left(v_0^{n_0} v_1^{n_1}
\right)^m(-2\frac{D_0}{n_0}+ 2\frac{D_1}{n_1})  \cr
 d \mapsto \frac{1}{n_0} D_0.
\end{array}
\end{eqnarray}
Note that if we set $T=A_{00} A_{11}v_0^{n_0}v_1^{n_1}$ this implies that
$$
\hat{F}^A_{n_0,n_1}\circ \Psi^{-1}(t^mS)=T^m\hat{F}^A_{n_0,n_1}\circ
\Psi^{-1}(S)\quad\forall S\in sl(2,{\bb C}).
$$
\end{ex}
\begin{ex}
If $\{H_1,H_2,X_{\pm i}:\,1\le i\le 3\}$ is a standard basis of $sl(3,{\bb
C})$ then
$\{d,t^kH_1,t^lH_2,t^mX_{i},t^nX_{-i}:k,l,m,n\in {\bb Z},\, 1\le i\le 3\}$
is a basis of $\tilde{L}\left(sl(3,{\bb C})\right)$. The Cartan matrix of
$\tilde{sl}(3,{\bb C})$ is
$$
\begin{pmatrix} 2&-1&-1\\ -1&2&-1\\-1&-1&2\\\end{pmatrix}
$$
and the solution matrices are of the form
$$
A=\begin{pmatrix}
A_{00}&\frac{1}{2}(-1-\varepsilon)A_{11}&\frac{1}{2}(-1+\varepsilon)A_{22}\\
\frac{1}{2}(-1+\varepsilon)A_{00}&A_{11}&\frac{1}{2}(-1-\varepsilon)A_{22}\\
\frac{1}{2}(-1-\varepsilon)A_{00}&\frac{1}{2}(-1+\varepsilon)
A_{11}&A_{22}\\\end{pmatrix}
$$
By \eqref{loopiso} and \ref{expembn}, the images of $\{d, H_1,H_2,X_{\pm
1},X_{\pm 2},t^{\pm 1}X_{\mp 3}\}$ under the embedding
$\hat{F}^A_{n_0,n_1,n_2}\circ
\Psi^{-1}:\tilde{L}\left(sl(3,{\bb C})\right)\rightarrow {\cal X}(M)$ are
given by the formulae:
\begin{eqnarray}
\label{eq:tildesl3}
&&\left\{ \begin{array}{l}
X_1 \mapsto A_{11} v_1^{n_1}
\left(\frac{1}{2}(-1-\varepsilon) \frac{D_0}{n_0}+\frac{D_1}{n_1} +
 \frac{1}{2}(-1+\varepsilon) \frac{D_2}{n_2}\right)  \cr
X_2 \mapsto A_{22} v_2^{n_2}
\left( \frac{1}{2}(-1+\varepsilon)
\frac{D_0}{n_0}+\frac{1}{2}(-1-\varepsilon) \frac{D_1}{n_1}+
\frac{D_2}{n_2}\right) \cr
\frac{1}{t}X_3\mapsto \frac{-1}{A_{00}v_0^{n_0}}
\left( \frac{D_0}{n_0} +\frac{1}{2}(-1-\varepsilon) \frac{D_1}{n_1}
+\frac{1}{2}(-1+\varepsilon)\frac{D_2}{n_2}\right)
\end{array}
\right.  \nonumber \\
 \nonumber \\
&&\left\{
\begin{array}{l}
X_{-1} \mapsto  \frac{-1}{A_{11}v_1^{n_1}}
\left(\frac{1}{2}(-1+\varepsilon) \frac{D_0}{n_0}+ \frac{D_1}{n_1}
+ \frac{1}{2}(-1-\varepsilon) \frac{D_2}{n_2}\right)
 \cr
X_{-2}\mapsto  \frac{-1}{A_{22}v_2^{n_2}}
\left( \frac{1}{2}(-1-\varepsilon)
\frac{D_0}{n_0}+\frac{1}{2}(-1+\varepsilon) \frac{D_1}{n_1}
 + \frac{D_2}{n_2}\right)\cr
tX_{-3}\mapsto A_{00}v_0^{n_0}
\left( \frac{D_0}{n_0}+\frac{1}{2}(-1+\varepsilon) \frac{D_1}{n_1}+
\frac{1}{2}(-1-\varepsilon)
\frac{D_2}{n_2} \right)
\end{array}
\right.
\nonumber \\
  \\
&&\left\{
\begin{array}{l}
d\mapsto \frac{1}{n_0}D_0\cr
H_1 \mapsto  -\frac{D_0}{n_0}+2\frac{D_1}{n_1}-\frac{D_2}{n_2} \cr
H_2 \mapsto  -\frac{D_0}{n_0}-\frac{D_1}{n_1}+2 \frac{D_2}{n_2}
\end{array}
\right.
\nonumber
\end{eqnarray}
Calculating the images of $X_{\pm 3}=\pm[X_{\pm 1},X_{\pm 2}]$  gives
\begin{align}
X_3 &\mapsto -A_{11}A_{22}v_1^{n_1} v_2^{n_2}\varepsilon
\left( \frac{D_0}{n_0} +\frac{1}{2}(-1-\varepsilon) \frac{D_1}{n_1}
+\frac{1}{2}(-1+\varepsilon)\frac{D_2}{n_2}\right)\\
X_{-3}&\mapsto \frac{\varepsilon}{A_{11}A_{22}v_1^{n_1} v_2^{n_2}}
\left( \frac{D_0}{n_0}+\frac{1}{2}(-1+\varepsilon) \frac{D_1}{n_1}+
\frac{1}{2}(-1-\varepsilon)
\frac{D_2}{n_2} \right)
\end{align}
This gives  the formulae for the restriction of
$\hat{F}^A_{n_0,n_1,n_2}\circ\Psi^{-1}:\tilde{L}\left(sl(3,{\bb
C})\right)\rightarrow{\cal X}(M)$ to the subset $sl(3,{\bb C})=\langle
H_1,H_2,X_{\pm i}:\,1\le i\le 3\rangle$. Note that
setting
$D_0=0$ one recovers the formulae \eqref{eq:sl3}. Now observe that if we set
$T=\varepsilon A_{00}A_{11}A_{22}v_0^{n_0}v_1^{n_1}v_2^{n_2}$ then
\begin{align}
\hat{F}^A_{n_0,n_1,n_2}\circ\Psi^{-1}(X_{-3})&=\frac{1}{T}
\hat{F}^A_{n_0,n_1,n_2
}\circ\Psi^{-1}(tX_{-3})\\
\hat{F}^A_{n_0,n_1,n_2}\circ\Psi^{-1}(X_{3})&=T\hat{F}^A_{n_0,n_1,n_2}
\circ\Psi^{-1}(\frac{1}{t}X_{3})\\
\hat{F}^A_{n_0,n_1,n_2}\circ\Psi^{-1}(S)(T)&=0\quad\forall S\in sl(3,{\bb C})\\
\hat{F}^A_{n_0,n_1,n_2}\circ\Psi^{-1}(d)(T)&=T.
\end{align}
From this it follows that
$\hat{F}^A_{n_0,n_1,n_2}\circ\Psi^{-1}:\tilde{L}\left(sl(3,{\bb
C})\right)\rightarrow{\cal X}(M)$ is given by the formulae for
$\hat{F}^A_{n_0,n_1,n_2}\circ\Psi^{-1}(d),\hat{F}^A_{n_0,n_1,n_2}\circ
\Psi^{-1}(H_{i}),\hat{F}^A_{n_0,n_1,n_2}\circ\Psi^{-1}(X_{\pm
i})$ above and
\begin{equation}
\hat{F}^A_{n_0,n_1,n_2}\circ\Psi^{-1}(t^mS)=T^m\hat{F}^A_{n_0,n_1,n_2}
\circ\Psi^{-1}(S)\,\forall S\in
sl(3,{\bb C}).
\end{equation}
\end{ex}


\begin{thebibliography}{99}

\bibitem[G-K]{GaKac} O. Gabber, V. Kac :  On defining relations of certain
infinite-dimensional Lie algebras,
{\it Bull. Amer. Math. Soc}, {\bf 5}, 185-189 (1981).

\bibitem[K]{Kac} V.Kac : `` Infinite dimensional Lie algebras'',
Cambridge University Press,
Third Edition, 1990.

\bibitem[Mac]{Mac} I.G. Macdonald:  Kac-Moody algebras , in ``Lie algebras
and related topics'',
R.V. Moody (ed), Conference Proceedings of the Canadian Mathematical
Society, Vol. 5,
American Mathematical Society, 1986.

\bibitem[Mo]{Mo} R.V. Moody :  A new class of Lie algebras, {\it Journal
of Algebra}, {\bf 10}, 211-230 (1968).

\bibitem[R]{R} M. Rausch de Traubenberg:   Fractional Supersymmetry and
Infinite Dimensional Lie Algebras, hep-th/0109106, 
Nucl. Phys. Proc. Suppl. {\bf 102} (2001) 256-263.


\bibitem[Se]{Se} J.-P. Serre :  `` Alg\`ebres de Lie semi-simples
complexes'', Benjamin, New York, 1966.


\end{thebibliography}
\end{document}